\documentclass[a4paper,11pt]{amsart}
\usepackage[latin1]{inputenc}
\usepackage{amsmath}
\usepackage{amsfonts}
\usepackage{amssymb}
\usepackage{hyperref}
\usepackage{multirow}
\usepackage{tabularx}

\newtheorem{proposition}{Proposition}

\catcode`\@=11

\catcode`\;=\active
\def;{\relax\ifhmode\ifdim\lastskip>\z@\unskip\fi\kern.2em\fi\string;\ \ignorespaces}

\catcode`\:=\active
\def:{\relax\ifhmode\ifdim\lastskip>\z@\unskip\fi\penalty\@M\ \fi\string:\ \ignorespaces}

\catcode`\!=\active
\def!{\relax\ifhmode\ifdim\lastskip>\z@\unskip\fi\kern.2em\fi\string!\ \ignorespaces}

\catcode`\?=\active
\def?{\relax\ifhmode\ifdim\lastskip>\z@\unskip\fi\kern.2em\fi\string?\ \ignorespaces}

\newif\ifguill  \guilltrue
\catcode`\"=\active
\def"{\ifguill``\guillfalse\else''\guilltrue\fi}

\catcode`\@=12
\catcode`\à=\active \defà{\`a}%
\catcode`\À=\active \defÀ{\`A}%
\catcode`\â=\active \defâ{\^a}%
\catcode`\é=\active \defé{\'e}%
\catcode`\è=\active \defè{\`e}%
\catcode`\ê=\active \defê{\^e}%
\catcode`\ë=\active \defë{\"e}%
\catcode`\ï=\active \defï{\"\i}%
\catcode`\î=\active \defî{\^\i}%
\catcode`\ô=\active \defô{\^o}%
\catcode`\ö=\active \defö{\"o}%
\catcode`\ù=\active \defù{\`u}%
\catcode`\û=\active \defû{\^u}%
\catcode`\ü=\active \defü{\"u}%
\catcode`\ç=\active \defç{\c c}%

\newcommand{\Q}{\mathbb{Q}}

\newcommand{\be}{\begin{enumerate}}
\newcommand{\ee}{\end{enumerate}}
\newcommand{\soe}{\geq} 
\newcommand{\ioe}{\leq} 

\newcommand{\vers}{\longrightarrow}

\newcommand{\pair}{\text{even}}
\newcommand{\impair}{\text{odd}}

\begin{document}

\title[The conjecture H]{A conjectured lower bound for the cohomological dimension of elliptic spaces
\\Recent results  in some simple cases.
\\ \today
}
\vskip 10mm
\author{Mohamed Rachid Hilali}
\address{D\'epartement de Math\'ematiques\\
         Faculté des sciences Ain Chok\\
         Universit\'e Hassan II,  Route d'El Jadida\\
         Casablanca\\
         Maroc}
\email{rhilali@hotmail.com}

\author{My Ismail Mamouni}
\address{Classes préparatoires aux grandes écoles d'ingénieurs\\
         Lycée Med V\\
         Avenue 2 Mars\\
		Casablanca\\
        Maroc}
\email{myismail1@menara.ma}

\vskip 20mm
\begin{abstract}
Here we prove some special cases of the following conjecture: that the sum of the Betti numbers of a 1-connected elliptic space is greater than the total rank of its homotopy groups. Our main tool is Sullivan's minimal model.
\end{abstract}



\maketitle
\section{Introduction} The subject of our research is the following result:

\textbf{Conjecture H (Topological version)}. \textit{If $X$ is a 1-connected elliptic space then }$\dim H^*(X,\Q)\soe \dim \left(\pi_*(X)\otimes\Q\right)$.

\vskip 3mm
For elliptic spaces $X$, it seems roughly that $\dim H^*(X,\Q)$ becomes larger when $\dim\pi_*(X)\otimes\Q$ does. That was the key idea for the first author to conjecture a lower bound for the cohomolgical dimension of an elliptic space in terms of the total rank of its homotopy groups. Let us recall  that for any 1-connected space $X$ of finite type, i.e., $\dim H^*(X,\Q)<\infty$ for all $k\soe 0$, there exists a commutative differential graded algebra $(\Lambda V,d)$, called \textit{Sullivan model} of $X$, which algebraically models the rational homoptopy type of the space, more precisely $H^*(\Lambda V,d)\cong H^*(X,\Q)$ as algebras, and $V\cong\pi_*(X)\otimes\Q$ as vector spaces. By $\Lambda V$ we mean the free commutative algebra generated by the graded vector space $V$, i.e., $\Lambda V=TV\diagup \langle v\otimes w-(-1)^{|v||w|}w\otimes v\rangle$, where $TV$ denotes the tensor algebra over $V$. $\Lambda ^nV$ denotes the set of elements of $\Lambda V$ of wordlength $n$ and $\Lambda ^{\soe n}V:=\bigoplus_{k\soe n}\Lambda ^kV$ denotes that  of elements of $\Lambda V$ of wordlength at least $n$. The differential $d$ of any element of $V$ is a "polynomial" in $\Lambda V$ with no linear term, i.e., $dV\subset\Lambda^{\soe 2}V$, we say that the model is \textit{minimal}. A space $X$ and its minimal Sullivan model are called \textit{elliptic} if both $V$ and $H^*(\Lambda V,d)$ are finite dimensional spaces. Because of this contravariant correspondance between spaces and their minimal models, the topological version of our conjecture admits the following algebraic interpretation:

\vskip 3mm
\textbf{Conjecture H (Algebraic version)}. \textit{If $\Lambda V$ is a 1-connected elliptic Sullivan minimal model then }$\dim H^*(\Lambda V,d)\soe \dim V$.

\vskip 3mm
We assume that the minimal model is simply connected, i.e., that the vector space $V$ has no generators in degrees lower than 2. This assumption is necessary in order to translate our algebraic results into topological ones, although it is not necessary for the algebraic results themselves. For more details about minimal Sullivan models of spaces we refer the reader to \cite{FHT}, (138-160).

\vskip 3mm
Our main motivation to believe that the conjecture H holds in the elliptic case, is that it is for pure spaces \cite{Hi} and for H-spaces, symplectic manifolds,.... \cite{HM}. Homogenous spaces are pure. Topological groups are H-spaces, and Kahler manifolds are symplectic.
\section{Results and proofs}
\begin{proposition}
If $X$ is a nilmanifold with $(\Lambda V,d)$ as a model, then $\dim H^*(\Lambda V,d)\soe \dim V$.
\end{proposition}
\textit{Proof}. We know from \cite{Dix} that a nilmanifold $X$ of dimension $n$ has $b_i\soe 2$ for $1\ioe i\ioe n-1$ and hence $\dim H^*(\Lambda V,d)\soe 2fd(X)$, where $fd(X)$ called the \textit{formal dimension} of $X$, denotes the largest integer $k$ such that $H^k(X,\Q)\neq 0$. We know also from (\cite{FH}-Proposition 2.6) that space of "type F" checks the inequality $fd(X)\soe \dim V$. Nilmanifolds are  of "type F", since they are $K(G,1)$ where $G$ is a nilpotent group.\hfill $\Box$

\vskip 3mm
\begin{proposition}
If   an elliptic minimal model $(\Lambda V,d)$ has an homogeneous-length differential and whose rational Hurewicz homorphism is non-zero in some odd degree. Then $\dim H^*(\Lambda V,d)\soe \dim V$.
\end{proposition}
\textit{Proof}. We say that $\Lambda V$ has differential of homegeneous-length $l$ if $d:V\vers \Lambda^l V$. We know from \cite{Lu} that under hypotheses above we have $\dim H^*(\Lambda V,d)\soe 2{\rm cat}_0(\Lambda V)$ and from \cite{FH2} that $\dim V^\pair\ioe \dim V^\impair\ioe {\rm cat}_0(\Lambda V)$. We can then conclude that $\dim V\ioe 2{\rm cat}_0(\Lambda V)\ioe \dim H^*(\Lambda V,d)$. \hfill $\Box$

\vskip 3mm
\begin{proposition}\label{coformal}
If  an elliptic  minimal model $(\Lambda V,d)$ has a differential, homogenous of length at least 3, then $\dim H^*(\Lambda V,d)\soe \dim V$.
\end{proposition}
\textit{Proof}. The cohomology of such spaces admits a second grading $H^*(\Lambda V,d)=\bigoplus\limits_{k\soe 1}H^*_k(\Lambda V,d)$, given by length of representative cocycle. We know from (\cite{Lu}-Theorem 2.2) that $H_k^*(\Lambda V,d)\neq 0$ for each $k=0,\cdots,e$ where $e=\dim V^\impair+(l-2)\dim V^\pair$. Then $\dim H^*(\Lambda V,d)\soe e\soe \dim V$, when $l\soe 3$.

\vskip 3mm
\begin{proposition}
If  an elliptic  minimal model $(\Lambda V,d)$ has a differnetial, homogenous of length 2 (i.e: coformal)  with odd degree generators only, i.e., $V^\pair=0$, then $\dim H^*(\Lambda V,d)\soe \dim V$.
\end{proposition}
\textit{Proof}. The proof is similar to that of that Proposition  \ref{coformal}

\vskip 3mm
\begin{proposition} If $X$ is an elliptic space wich has the homotopy type of the $r$-product of elliptic spaces satisfying the conjecture H, then it is also for $X$.
\end{proposition}
\textit{Proof}. The argument is that: $\dim H^*(Y\times Z,d)= \dim H^*(Y,d).\dim H^*(Z,d)$ and that $\dim \left(\pi(Y\times Z)\otimes\Q\right)=\dim \left(\pi(Y)\otimes\Q\right)+\dim \left(\pi(Z)\otimes\Q\right)$. \hfill $\Box$

\vskip 3mm
\begin{proposition} If $(\Lambda V,d)$ is a formal and hyperelliptic model, then $\dim H^*(\Lambda V,d)\soe \dim V$.
\end{proposition}
\textit{Proof}. The Sullivan model $(\Lambda V,d)$ is called hyperelliptic if $dV^\pair=0$ and $dV^\impair\subset \Lambda^+V^\pair\otimes \Lambda V^\impair$. Set $\Lambda V^\pair=\Lambda \{x_1,\ldots,x_n\}$ and $V^\impair=\Lambda \{y_1,\ldots,y_{n+p}\}$ with $p=\dim V^\impair-\dim V^\pair\soe 0$ (cf. \cite{FHT}-Proposition 32.10 (444)). The first author  (cf. \cite{Hi}) has showed that the conjecture H  holds for pure spaces,  then  we assume that $(\Lambda V,d)$  is not pure. As it is formal,  then $\dim H^*(\Lambda V,d)\soe 2^p$ when $H^*(\Lambda V,d)$ is given by $\Q\left[x_1,\ldots,x_n \right] \otimes \Lambda \{y_1,\ldots,y_{n}\}\diagup (dy_{n+1},\ldots,dy_{n+p})$ with $dx_i=dy_i=0$ for $1\ioe i\ioe n$. On the other hand it is proved (cf. \cite{HM}-Theorem C) for hyperelliptic models that $\dim H^*(\Lambda V,d)\soe \dim V$ when $\dim H^*(\Lambda V,d)\soe 2^{p-1}$, this simple remark completes the proof. \hfill $\Box$

\vskip 3mm
\begin{proposition} If  $(\Lambda V=\Lambda (U,W),d)$ is a two-stage, elliptic minimal model with odd degree generators only and suppose that $d: W\vers \Lambda ^2U$ is an isomorphism, then $\dim H^*  (\Lambda V,d)   \soe \dim V$.
\end{proposition}
\textit{Proof.} A minimal model $(\Lambda V,d)$ is said to be \textit{two-stage} if $V$ decomposes $V=U\oplus W$ with $dU=0$ and $dW\subset \Lambda U$, then the Mapping theorem (\cite{FHT}, page 375) implies that $W$ has generators of odd degree only. On the other hand, $U$ may have generators of odd or even degree. We know from (\cite{JL}-Proposition 2.1), that under hypotheses here  above we have $\dim H^*  (\Lambda (U,W),d)   \soe 2^{\dim W}$. Set $\dim U=n$, then $\dim W=\dfrac{n(n+1)}2\soe n=\dim U$ and $\dim W\soe \dfrac{\dim V}2=m$. Finally $\dim H^*  (\Lambda V,d) \soe 2^m\soe 2m=\dim V$. \hfill $\Box$

\vskip 3mm
\textbf{\underline{Open Questions:} } 
\vskip 4mm
\begin{itemize}
\item The author in \cite{Lu} believes that his Conjecture 3.4 holds at least in the coformal case. Has this been settled?
\item We know from (\cite{JL}-Corollary 3.5) that if $X$ has  a two-stage model with odd degree generators only, then $\dim H^*(X,\Q)\soe 2^{\dim G_*(X)}$, where $G_*(X)$ denoted the subgroup of $\pi_*(X)$ called the \textit{Gottlieb group}. A question is: are there any cases where  the equality $\dim G_*(X)= \dim \left(\pi_*(X)\otimes\Q\right)$ holds ?
\item If $F\vers E\vers B$ is a fibration where $F$ and $B$ are elliptic and both verify the conjecture H, what conditions on the fibration will guarantee that $E$ will too?
\end{itemize}

\vskip 4mm
\textbf{Acknowledgements}. It is for us a pleasure to thank Micheline Vigué (Univ. Paris 13, French) and Barry Jessup (Univ. Ottawa, Canada)  for their  interest and for their several readings and corrections.

\vskip 4mm
\end{document}